\documentclass[12pt]{article}
\usepackage{geometry}                
\usepackage[parfill]{parskip}    
\usepackage{graphicx}
\usepackage{amsmath, latexsym, amssymb}
\input xypic
\begin{document}
\title{Yang-Mills bar connections over compact K\"ahler manifolds}
\author{\bf  H\^ong V\^an L\^e\\
Mathematical Institute of ASCR,\\
Zitna 25, CZ-11567 Praha 1, \\
Czech republic\\
email: hvle@math.cas.cz}

\date{}
\maketitle
\newcommand{\R}{{\mathbb R}}
\newcommand{\C}{{\mathbb C}}
\newcommand{\F}{{\mathbb F}}
\newcommand{\Z}{{\mathbb Z}}
\newcommand{\N}{{\mathbb N}}
\newcommand{\Q}{{\mathbb Q}}

\newcommand{\Aa}{{\mathcal A}}
\newcommand{\Bb}{{\mathcal B}}
\newcommand{\Cc}{{\mathcal C}}    
\newcommand{\Dd}{{\mathcal D}}
\newcommand{\Ee}{{\mathcal E}}
\newcommand{\Ff}{{\mathcal F}}
\newcommand{\Gg}{{\mathcal G}}    
\newcommand{\Hh}{{\mathcal H}}
\newcommand{\Kk}{{\mathcal K}}
\newcommand{\Ii}{{\mathcal I}}
\newcommand{\Jj}{{\mathcal J}}
\newcommand{\Ll}{{\mathcal L}}    
\newcommand{\Mm}{{\mathcal M}}    
\newcommand{\Nn}{{\mathcal N}}
\newcommand{\Oo}{{\mathcal O}}
\newcommand{\Pp}{{\mathcal P}}
\newcommand{\Qq}{{\mathcal Q}}
\newcommand{\Rr}{{\mathcal R}}
\newcommand{\Ss}{{\mathcal S}}
\newcommand{\Tt}{{\mathcal T}}
\newcommand{\Uu}{{\mathcal U}}
\newcommand{\Vv}{{\mathcal V}}
\newcommand{\Ww}{{\mathcal W}}
\newcommand{\Xx}{{\mathcal X}}
\newcommand{\Yy}{{\mathcal Y}}
\newcommand{\Zz}{{\mathcal Z}}

\newcommand{\zt}{{\tilde z}}
\newcommand{\xt}{{\tilde x}}
\newcommand{\Ht}{\widetilde{H}}
\newcommand{\ut}{{\tilde u}}
\newcommand{\Mt}{{\widetilde M}}
\newcommand{\Llt}{{\widetilde{\mathcal L}}}
\newcommand{\yt}{{\tilde y}}
\newcommand{\vt}{{\tilde v}}
\newcommand{\Ppt}{{\widetilde{\mathcal P}}}
\newcommand{\bp }{{\bar \partial}} 

\newcommand{\Remark}{{\it Remark}}
\newcommand{\Proof}{{\it Proof}}
\newcommand{\ad}{{\rm ad}}
\newcommand{\Om}{{\Omega}}
\newcommand{\om}{{\omega}}
\newcommand{\eps}{{\varepsilon}}
\newcommand{\Di}{{\rm Diff}}
\newcommand{\Pro}[1]{\noindent {\bf Proposition #1}}
\newcommand{\Thm}[1]{\noindent {\bf Theorem #1}}
\newcommand{\Lem}[1]{\noindent {\bf Lemma #1 }}
\newcommand{\An}[1]{\noindent {\bf Anmerkung #1}}
\newcommand{\Kor}[1]{\noindent {\bf Korollar #1}}
\newcommand{\Satz}[1]{\noindent {\bf Satz #1}}

\newcommand{\gl}{{\frak gl}}
\renewcommand{\o}{{\frak o}}
\newcommand{\so}{{\frak so}}
\renewcommand{\u}{{\frak u}}
\newcommand{\su}{{\frak su}}
\newcommand{\ssl}{{\frak sl}}
\newcommand{\ssp}{{\frak sp}}

\newcommand{\Cinf}{C^{\infty}}
\newcommand{\CS}{{\mathcal{CS}}}
\newcommand{\YM}{{\mathcal{YM}}}
\newcommand{\Jreg}{{\mathcal J}_{\rm reg}}
\newcommand{\Hreg}{{\mathcal H}_{\rm reg}}
\newcommand{\SP}{{\rm SP}}
\newcommand{\im}{{\rm im}}

\newcommand{\inner}[2]{\langle #1, #2\rangle}    
\newcommand{\Inner}[2]{#1\cdot#2}
\def\NABLA#1{{\mathop{\nabla\kern-.5ex\lower1ex\hbox{$#1$}}}}
\def\Nabla#1{\nabla\kern-.5ex{}_#1}

\newcommand{\half}{\scriptstyle\frac{1}{2}}
\newcommand{\p}{{\partial}}
\newcommand{\notsub}{\not\subset}
\newcommand{\iI}{{I}}               
\newcommand{\bI}{{\partial I}}      
\newcommand{\LRA}{\Longrightarrow}
\newcommand{\LLA}{\Longleftarrow}
\newcommand{\lra}{\longrightarrow}
\newcommand{\LLR}{\Longleftrightarrow}
\newcommand{\lla}{\longleftarrow}
\newcommand{\INTO}{\hookrightarrow}

\newcommand{\Sy}{\text{ Diff }_{\om}}
\newcommand{\Ex}{\text{Diff }_{ex}}
\newcommand{\jdef}[1]{{\bf #1}}
\newcommand{\QED}{\hfill$\Box$\medskip}

\newcommand{\UuU}{\Upsilon _{\delta}(H_0) \times \Uu _{\delta} (J_0)}
\newcommand{\bm}{\boldmath}
\medskip


\abstract  In this note we introduce a Yang-Mills  bar  equation on complex vector  bundles  $E$ provided
with a  Hemitian  metric
over compact Hermitian manifolds.
According to the  Koszul-Malgrange criterion any holomorphic structure on $E$ can be seen as
a solution to this equation. We show  the  existence  of a non-trivial solution to this equation over compact 
K\"ahler manifolds as well as a short time existence of a related negative Yang-Mills bar gradient flow. We   also show  a rigidity of holomorphic connections among a class  of  Yang-Mills bar connections over  compact  K\"ahler   manifolds of positive Ricci curvature.
\endabstract



\section {Introduction}

Let $M^{2n}$ be a compact Hermitian manifold of  real dimension $2n$ and $E$ be a complex vector bundle  over $M ^{2n}$.
 The following  Koszul-Malgrange  criterion \cite{K-M1958}, see also \cite{D-K1990}, 2.1.53, 2.1.54,  establishes  the equivalence
between the  existence of a holomorphic structure on  $E$  and   a partial flatness of   $E$.

\medskip

{\bf Koszul-Malgrange criterion.}  {\it A  complex
vector bundle $E$ over a complex manifold $M^{2n}$ carries a holomorphic structure, if and only if there
is a  connection $A$ on $E$  such that the $(0,2)$-component  $F ^{0,2}_A$ of
the curvature $F_A$ of $A$ vanishes.}

\medskip

Thus we shall call  a  connection $A$ satisfying the Koszul-Malgrange criterion  a  holomorphic connection.
It is well-known (see e.g. \cite{D-K1990}) that we can  replace the  connection $A$ in the  Kozsul-Malgrange criterion
by a unitary connection $A$  for any given choice of a  compatible (Hermitian) metric $h$ on $E$.
\medskip

We  introduce in  section 2 (see (2.5.1) and (2.5.2))  a Yang-Mills bar equation as the Euler-Lagrange equation for
the Yang-Mills bar functional  which is   the square of the $L_2$-norm of the $(0,2)$-component
$F^{0,2}_A$ of  a   unitary  connection $A$ on $(E,h)$. Solutions of a Yang-Mills bar equation are called Yang-Mills bar connections. 
The Yang-Mills bar equation has an advantage  over the equation for a holomorphic connection,   because the later one is overdetermined if the complex dimension of the bundle is greater or equal to 2 and  $n \ge 4$, and
the  first one is elliptic modulo a degeneracy  which is  formally  generated by an action of the   complex gauge group of the  complex vector bundle $E$ (the degeneracy is   formal generated since the action of this group on the   ``small" space does not  preserve the Yang-Mills bar functional,  see   2.7.b and Remark 5.13). Thus we hope  that  by using this equation we  shall be able 
to  find  useful sufficient conditions  under which a complex vector bundle carries a holomorphic structure.  Appropriate sufficient conditions  for the existence of a holomorphic structure on complex vector bundles over
 projective algebraic manifolds  could  be a key step in solving the Hodge conjecture, if the conjecture is correct. A particular  result in this direction is  our Theorem  4.25 which states that an almost holomorphic   connection  over a compact K\"ahler manifold of positive Ricci curvature  is holomorphic, in particular any Yang-Mills bar connection
 on a  4-dimensional compact K\"ahler manifold of positive Ricci curvature is holomorphic.  

In section 2 after introducing the Yang-Mills  equation we also discuss the symmetry  of this equation in 2.7. In section 3 we  give a proof of  the Hodge-K\"ahler identities
for general  unitary connections over   K\"ahler manifolds and show the existence of
non-trivial  Yang-Mills bar connections. In section 4 we derive a Bochner-Weitzenb\"ock
type identity on compact K\"ahler manifolds and  prove  Theorem 4.25. In section
5 we  introduce the notion of affine integrability condition,  a negative Yang-Mills bar
gradient flow and  find  an affine integrability condition for this flow (Theorem 5.9).
Unlike   previously known cases for  weakly parabolic equations (Ricci flow, Yang-Mills flow),  our  affine integrability  is not derived
from  an action of a group, which preserves the Lagrangian on the  space, where our flow is considered (see 2.7.b and Remark 5.13.i). The automorphism group of the  Yang-Mill bar equation gives us only "half" of the integrability condition.  In particular, the  DeTurck approach to weakly parabolic equations  seems    inapplicable  to our flow. 
 In the last section 6 we  prove the short time existence, uniqueness and smoothness of  a solution of
an evolution equation with affine integrability condition, slightly  extending  a Hamilton's result.

 \medskip

\section { Yang-Mills bar equation}

Let $(V, \langle,\rangle)$ be a Euclidean space. Denote by $V_\C$ its complexification.  Then $\langle,\rangle$ extends uniquely  to a complex bilinear   form $\langle,\rangle_\C : V_{ \C}
\times V_{ \C} \to  \C$.  Denote by $(v,w) : =   \langle v,\bar w\rangle_\C $ the  associated Hermitian form on $V _\C$ and  by
$\langle v, w\rangle = Re  ( v,w)$  the Euclidean metric on the
space $(V_{ \C})\otimes \R$. We note that the restriction  of this metric to $V$ coincides with the original metric
$\langle,\rangle$. 
Conversely any Hermitian metric ($J$-invariant  Euclidean metric)
on  a complex space $(V, J)$ considered as a complexification of a real vector space $V _0$ is obtained in this way.

In this note  we shall define by the same $(,)$ (and resp. $\langle,\rangle$) the  Hermitian form (resp. the Euclidean metric) extended in the above way
 from any   vector bundle $(E, \langle,\rangle)$ provided  with a fiber-wise Euclidean metric $\langle,\rangle$ to its complexification $E_{ \C}$ (resp. considered as a  real space).
 If $A$ is a connection on  $(E, \langle,\rangle)$ then $A$  can be extended to a  unitary connection also denoted by $A$ on the
 complexification $E_\C$ 
  with that extended metric by setting $d _A  ( \sqrt{-1} \phi) : = \sqrt{-1}  d_A  (\phi)$.
 
 \medskip
 
 Now let $A$ be a connection on a complex vector bundle $(E,J)$ over a Hermitian manifold $M^{2n}$. Denote by $\Om^{ p,q}  (E)$ the space
 of $E$-valued $(p,q)$-forms  on $M^{2n}$: $\Om^{p,q} ( E) = \Om^{p,q} (M) \otimes _\C E$. 
  We have  the decomposition
 $$ d_A = \p _A \oplus \bar \p _A : \Om (E) \to \Om ^ {1,0} (E) \oplus \Om ^{0,1} (E).$$
   In general we  have the inclusion 
 $$d _A ( \Om ^{ p,q} ( E) ) \subset \Om ^{ p+1, q} ( E) \oplus \Om ^{ p, q +1} (E),$$
 since for $\psi \in  \Om ^0 ( E)$ and $\phi \in \Om ^{p,q}(M^{2n}) $ we have
 $$d _A ( \psi \otimes \phi) = d_A  (\psi) \otimes \phi + \psi\otimes d\phi  \in \Om ^{ p+1, q} (E) \oplus \Om ^{p, q+1}(E).$$
(The operator $d_A$ is well defined on  $\Om ^{p,q} (E)$, since  $d_A (J \psi)  = J d_A (\psi)$.)  For $\phi \in \Om ^{ p,q}(E)$ we shall denote by $ \p_A(\phi)$  the projection of  $d_A (\phi)$ on the first
 factor and by $\bar \p_A  (\phi)$ the projection on the second factor w.r.t the above  decomposition.

 We note that the curvature  $F_A\in \Om ^2 (End_J E)$ of $A$   can be considered as an element in $\Om ^2 _{\C}  ( End_J (E))$ 
 
 Let $(E,h)$ be a Hermitian vector bundle, i.e. a complex vector bundle $(E,J)$ provided with  a Hermitian metric  $h$ but $E$
 need not to be holomorphic.  There is a natural (Killing) metric on the space
 $u_E$  of  skew-Hermitian endomorphisms of $E$, defined by  $\langle \theta_1, \theta_2 \rangle  =   - Re\, Tr (\theta_1 \cdot \theta_2)$.  We can also write
 $End_J E = u _E \oplus \sqrt{-1} u_E$. Thus  the metric $h$  extends to a positive definite
 bilinear form on $End_JE$ (defined by $\langle \theta_1, \theta_2 \rangle = Re \, Tr (\theta_1\cdot \theta_2^*)$). Here  $ \theta^*$ is the  conjugate   transpose of $\theta$, the adjoint of $\theta$ w.r.t the unitary
 metric $h$.   We note  that
 this metric is invariant under the  original   complex structure on $End _J (E)$ induced by $J$ which
 we denoted above by multiplication with $\sqrt{-1}$. Hence by the remark at the beginning of the section,  this metric extends
  to a metric on the space $\Om ^k_\C ( End_J E)$  by combining  the Killing  metric with
  the Hermitian metric on $M^{2n}$. The decomposition $\Om^k_\C (End _J E) = \sum _{ p + q = k} \Om^{p,q} (End_J E)$ is  an orthogonal decomposition w.r.t this metric.

 If $A$ is a  unitary connection on $(E,h)$, then  $F_A \in  \Om ^2 ( u_E) \subset \Om^2 (End_J E)$.
 We also note that in   the decomposition for the curvature of  unitary connection $A$:
 $$ F_A = (F_A)^{2,0}  + (F_A)^{1,1} + (F_A) ^{ 0,2} $$
 we have $ (F_A) ^{0,2} = - ( (F_A)^{2,0} ) ^ *$.  The Kozsul-Malgrange  criterion  suggests us to consider the following  Yang-Mills bar functional  on the space of all  unitary   connections $A$ on $(E,h)$ over $M^{2n}$
$$\Yy\Mm ^ b (A) = (1/2) \int_{M^{2n} } || (F_A) ^{0,2} ||.$$

It is easy to see that the functional $\Yy\Mm^b$ is invariant under the gauge transformation of the Hermitian  vector bundle $(E,h)$. 
We shall  derive the first variation formula for the Yang-Mills bar equation.  First we shall extend the  usual Hodge operator
 $*:\Om ^p (M^{2n}) \to \Om ^{2n-p} (M^{2n})$ to $\bar * : \Om ^{p}(End_J E) \to \Om ^{2n-p}(End_J E)$ defined as follows. 
 We extend $\bar *: \Om ^{p}(End_J E) \to \Om ^{2n-p} (End_J E)$ so that for each
 $\alpha \in \Om ^p_\C (End_J E) $ and $\beta 
 \in \Om ^p_\C (End_J E)$ we have
 $$ \langle \alpha(x),  \beta(x) \rangle = \langle vol_x M^{2n} , \alpha(x) \wedge^{(,)} (\bar * \beta(x))\rangle.\leqno (2.1)$$
 Here $\wedge ^{(,)}$ denotes  the composition of  the wedge product  with the contraction of the coefficients in  $End_J E$   via  the natural Hermitian form
 $(, )$ on $End_J E$. 
 
 Next we note  that  $A$ induces  naturally a  connection,  also denoted by $A$,  on the  Hermitian  vector
 bundle $End_JE$ provided with  the metric described above. It is known that  the curvature  $F_A$  of this induced  connection acts on the
 space $\Om ^0 (End _J E)$ as follows
 $$ F_A (\phi) = F _A \wedge \phi  : = [ F_A, \phi], \leqno (2.2)$$
 see e.g. \cite{B-L1981}, (2.7). (The wedge product  of  differential forms with coefficients in a Lie algebra
 bundle is  the composition of  the wedge product and   the Lie bracket).
 
Now we define the operator $\bp_A ^* : \Om ^{p,q} (End_J E) \to \Om ^{p, q-1}(End_J E)$ as follows (see also \cite{Kobayashi1987}, chapter III, (2.19), or \cite{G-H1994}, chapter 1, \S 2, for
the case that $E$  is absent)
  $$ (\bp _A ^*)\beta ^{p,q}: = (-1)  \bar * \bp_A \bar * \beta ^{ p,q}.\leqno (2.3)$$
Using the  following identity  for the formal adjoint $d_A ^ *$ of $d _A$  on  an even dimensional manifold $M^{2n}$  (see e.g. \cite{B-L1981}, (2.27), for the real case,  the complex case can be proved by the same way by  
using the Stocks formula  locally):
$$ (d_A ^ * ) \beta =  ( -1) \bar * d_A \bar *\beta $$
and taking into account (2.3) which implies that $\bp _A  ^ *$ is the component with  correct bi-degree  of  $d  _A ^ *$, we  conclude that  $\bp_A ^ *$ is the formal  adjoint of $\bp _A $. 
Now using the formula $(F_{A + ta}) ^{0,2} = (F_A)^{0,2} + t\bar \p _A a^{ 0,1} +  t ^2 a ^{0,1} \wedge a^{0,1}$ and taking into account  (2.2)  we get  immediately

\medskip
{\bf  2.4. Lemma.} {\it Let $M^{2n}$ be a compact Hermitian manifold  with (possibly empty) boundary. The first variation of the  Yang-Mills bar functional  is given by the formula}
$${d\over dt}_{| t = 0}   \Yy\Mm ^b(A+ ta)  =  \int_{M^{2n}} \langle(\bar \p_A ) ^* F^{0,2}_A, a \rangle   + \int_{\p M^{2n}} \langle vol_x,  a \wedge  ^{ (,)}\bar * F ^{0,2}_A\rangle.$$

\medskip

We shall call a smooth unitary connection $A$  a Yang-Mills bar connection, if   it satisfies the following two conditions
$$(\bar \p _A)^* F^{0,2}_A = 0,\leqno (2.5.1)$$
$$( \bar * F^{0,2}_A) _{|\p M^{2n}}= 0.\leqno (2.5.2)$$

Let $\triangle _A ^{\bar \p} : = \bar \p_A ( \bar \p_A)  ^* + (\bar \p_A)^* \bar \p_A$.
Using the Bianchi identity $\bar\p _A F ^{ 0,2}_A = 0$, which follows from the usual Bianchi identity,  and  using  the  equality $\langle vol_x, a \wedge ^{(,)}\bar * b \rangle = \langle vol_x, b \wedge^{(,)} \bar *a \rangle$,  we conclude that  we can replace  (2.5.1) in the system of two equations (2.5.1) and (2.5.2) by the following condition
$$\triangle ^{\bar \p} _A (F_A)^{0,2} = 0,\leqno (2.6.1)$$
to get an equivalent system of equations.

\medskip

{\bf 2.7. Symmetries of the Yang-Mills bar equation. }  a) We can   vary the  Yang-Mills bar functional among
all compatible  Hermitian metrics $h'$ on $(E,J)$  in order to get an  invariant of  the  complex vector bundle $E$.   Let $A_t$ be  a family of unitary connections w.r.t. a compatible metric
$h_t$. We note that we can write $ h _t = g_t (h)$, where $g_t$ is a (complex) gauge transformation of $(E, J)$. Clearly $ (g_t) ^{-1} A_t$ is a unitary connection w.r.t. $h$ ( i.e. $d_{( g_t ) ^{-1} A_t} h = 0 $).
 Now we have $ F  ^{0,2} _{ A_t} =  Ad_{ g_t}F  ^{ 0,2} _{ (g_t) ^{-1} A_t }$.  Moreover
 $$ || F ^{0,2} _{ A_t} || _{ h_t } = || Ad_{ g_t} ^{-1}  F^{ 0,2}_{ A_t } || _h  = || F ^{0,2} _{(g_t)^{-1}( A_t) } || _h.\leqno (2.7.1)$$
 (We can get (2.7.1) easily by noticing that the inner products on $End_JE$ induced
 by $h$ and $g (h)$
 satisfy the following relation
 $$\langle A, B \rangle _ { g ( h)} =  \sum _i \langle A ( g (e_i)), B(g(e_i)) \rangle _{ g (h)}  = \sum _i \langle Ad_{ g^{-1} } A ( e_i), Ad _{ g^{-1} } B ( e_i) \rangle _h$$
 where $e_i$ is an orthonormal basis in  $E$ w.r.t. $h$.)
  
 Hence the infimum of the Yang-Mills bar functional   is a constant which does not depend on the unitary metric $h$.  
 
 b)  The linearization of  the Yang-Mills bar equation is not elliptic  because  the  equation is invariant under the    gauge group $\Gg(E,h)$ of $(E,h)$, see (2.7.1).  The complexification of  this group is the gauge group $\Gg (E)$. This  complexified group acts  also on the space  $\Aa (E,h)$ of  all unitary connections w.r.t. a fixed  compatible metric $h$  \cite{D-K1990}, (6.1.4). For $g \in \Gg (E)$  we denote by $\hat g$  the new  (non-canonical) action of $g$  
 on  $\Aa (E,h)$ defined as follows
 
$$\bar \p _{ \hat g ( A) }  = g \bar \p _A g ^{-1} =  \bar \p_A - ( \bar \p _A  g) g ^{-1}, $$
 $$\p _{ \hat g ( A) }  =   \p_A +[ ( \bar \p _A  g) g ^{-1}] ^*. $$
 
 Though this action of $\hat \Gg (E)$ does not preserve the Yang-Mills bar functional, infinitesimally   it fails to do it at a connection $A$ only by a quadratic term in $F^{0,2}_A$ (see (5.3)).    
 \medskip

\section{  Yang-Mills bar connections over compact  K\"ahler manifolds}

Suppose that $A$ is a  unitary connection on a Hermitian vector bundle $E$ over
a  K\"ahler manifold $M^{2n}$   with a K\"ahler form $\om$. As before denote by $\bar \p _A ^*$ the formal adjoint of
$\bar \p _A : \Om^{p,q} (E) \to \Om^{p,q+1}  (E)$ defined by (2.3),  and by $\p _A ^*$ the formal
adjoint  of $\p _A: \Om^{p,q}   (E) \to  \Om^{p+1, q}  (E)$ defined in the same way.

Denote by $\Lambda: \Om^{p,q} (E) \to \Om ^{p-1, q-1}  (E)$ the adjoint
of the wedge multiplication by $\om$, an algebraic operator.
The following Hodge-K\"ahler identities
$$\bar \p _ A ^*= \sqrt{-1} [\p_A, \Lambda], \leqno (3.1)$$
$$\p _A^* = - \sqrt{-1} [\bar \p _A , \Lambda], \leqno (3.2)$$
are well-known for    the case of a holomorphic bundle $E$  and $A$ being   a unitary holomorphic connection
\cite{Kobayashi1987}, chapter III, (2.39), or  \cite{G-H1994}, chapter 0, \S 7, chapter 1, \S 2, where they are called the Hodge identities. These identities have been called K\"ahler identities in  \cite{D-K1990},  \S 6.1.   We shall  give a proof of the general case here, assuming  the validity of  (3.1) and (3.2)
 for unitary holomorphic connections $A$.
 
 Note  that
 it suffices to prove these identities locally, so we can assume that the bundle is $U(n)$-trivial and
 $\p_A = \p + A^{1,0}$, where $A ^{1,0} =\sum_{ i = 1} ^n A _i d z_i$, $ A_i\in End _J (E)$.
 Similarly
 $\bar \p_A= \bar \p +  A ^{0,1}$ with $A ^{1,0} =\sum_{i =1} ^n- ( A_i ) ^* d\bar z_i$.
 Here we define $\bar \p$ and $\p$ to be the $(1,0)$ and $(0,1)$ components of the unique unitary connection  which is compatible with the trivial holomorphic structure.  
 
 Since  the Hodge-K\"ahler identities  are valid for $A =0$,
it is easy to see that  (3.1) and (3.2)  are equivalent  to  
   following  algebraic identities
$$[ A^{0,1}] ^* = \sqrt{-1} [ A ^{1,0}, \Lambda], \leqno (3.3)$$
$$[ A^{1,0}] ^*  = -\sqrt{-1} [  A^{0,1}, \Lambda]. \leqno (3.4)$$
In view of the Hermitian linearity of  LHS  of (3.3) and (3.4):
$$ ( \lambda A + \gamma B)  ^* = \bar \lambda A ^* + \bar \gamma B ^*$$
 for $\lambda, \gamma \in \C $, and taking into account  the unitary of $A$ which implies $A ^{1,0} = - (A ^{ 0,1})^*$, it suffices to prove
 these identities for   a  $\C$-basic $\{ A^{1,0}  =e_{ij} dz_k,|\: 1 \le i, j \le \dim _\C E, \, 1\le k \le \dim_\C M^{2n} = n\} $  of  $(0,1)$-forms in
 $\Om^{0,1}  (End _J E)$. Here  $e_{ij}$ is an elementary
 matrix  in $End _J (E)$.
 We also assume that the K\"ahler metric at a given point $x$  is $\sum_i dz _i d\bar z_i$.
Denote by $i_k$ and $\bar i_k$ the adjoint of the multiplication
operators $dz_k \wedge$ and $d\bar z_k\wedge$ correspondingly. Then we have
$$[A^{1,0}] ^* =  (e_{ji}i_k), \:  [  A ^{0,1}]^* = -( e_{ij} \bar i_k)$$
$$\Lambda= -{\sqrt {-1}\over 2} \sum_{k=1} ^n \bar i_k i_k.$$
Substituting these identities in  LHS of (3.3) and (3.4)  we conclude that (3.3) and
(3.4) are  equivalent to the following identities for all $i,j,k$
$$- ( e_{ij} \bar i_k ) = \sqrt{-1} [e_{ij}dz_k, -{ \sqrt{-1}\over 2}\sum_{ k =1} ^n \bar i_k i_k] , \leqno (3.5)$$
$$ ( e_{ji} i_k) = -\sqrt{-1} [ - e_{ji}d\bar z_k, -{ \sqrt{-1}\over 2}\sum_{k =1} ^n \bar i_k i_k] , \leqno (3.6)$$
In their turn  (3.5) and (3.6)  are  immediate consequences of the following identities
$$ - \bar i_k = {1\over 2} [ d  z_k\wedge,  \sum_{ k =1} ^n \bar i_k i_k].\leqno (3.7)$$
$$-i_k = {1\over 2}[ d\bar z_k \wedge, \sum_{ k =1} ^n \bar i_k i_k], \leqno (3.8)$$
To prove (3.7) (and  (3.8)  resp.) we shall  compare the action of LHS of (3.7) (and of (3.8) resp.) and the action of RHS of (3.7) (and of (3.8) resp)  on $ \phi =  dz_J \wedge d\bar z_K$.
We use the following formulas proved in   p.112-113
of \cite{G-H1994}
$$ i_k ( dz_J \wedge d\bar z_K) = 0, \text { if } k \not \in J, \leqno (3.9)$$
$$i_k ( dz_k\wedge dz_J \wedge d\bar z_K) = 2 d z_J \wedge d\bar z_K,\leqno (3.10)$$
$$\bar i_k ( dz_J \wedge  d \bar z _K) = 0, \text { if } k \not \in K,\leqno (3.11)$$
$$\bar i_k (d\bar z_k \wedge  dz_J \wedge d\bar z_K) = 2 dz_J \wedge d\bar z_K.\leqno (3.12)$$
With help of (3.9) -(3.12) we get (3.7) immediately. It is easy to see
that (3.8) can be obtained from (3.7) by changing the complex orientation.
 \QED



 
 \medskip

 Set $\triangle _A ^\p : = \p _A \p_A^* + \p_A^* \p_A$, $\triangle _A ^{\bar \p} : = \bar \p _A \bar \p_A^* + \bar \p_A^* \bar \p_A$
 
 \medskip
 
{\bf 3.13. Corollaries}. {\it   For  $\phi, \psi \in \Om  ^{0,p } (  E )$ 
 we have the following simple expressions
 $$  \bar \p _A ^*\phi = -\sqrt{-1}  \Lambda \p _A(\phi),\leqno (3.13.1)$$ 
$$\int_{M^{2n}} \langle\sqrt{-1}\Lambda  F^{1,1} _A  \phi, \psi\rangle = \int_{M^{2n}}- \langle  \bar \p _A^*  \phi, \bar \p _A  ^* \psi\rangle +  \langle \p_A \phi, \p_A  \psi \rangle - \langle\bar \p _A \phi, \bar \p_A \psi\rangle .\leqno (3.13.2)$$
More generally, for all $\phi \in \Om ^{p,q} ( E)$ we have}
$$(\triangle _A^\p - \triangle _A ^{\bar \p})\phi= -\sqrt{-1} [ F^{1,1}_A \wedge,  \Lambda]\phi\leqno (3.13.3) $$
$$ \triangle_A ^{\bar \p}\phi = {1 \over 2} ( \triangle _A^d   + \sqrt{-1}[  -F^{0,2}_A + F ^{2,0}_A +F_A^{1,1}, \Lambda ])\phi\leqno (3.13.4)$$

\medskip

{\it  Proof.} 1) The first statement follows  immediately from the Hodge-K\"ahler identity  (3.1).
\medskip

 2) Substituting  $F^{1,1}_A = \bar \p _A \p _A + \p_A \bar \p _A$ we get
 $$\int_{ M^{2n} } \langle \sqrt{-1}\Lambda F ^{ 1,1} _A \phi, \psi \rangle = \int_{ M^{ 2n}} \langle\sqrt{-1}\Lambda ( \bar\p_A \p_A + \p_A\bar \p_A )\phi, \psi\rangle. $$
 
Now  applying the Hodge-K\"ahler identities to this  equation we get
 $$\int_{ M^{2n} } \langle\sqrt{-1} \Lambda F ^{ 1,1} _A \phi, \psi \rangle = \int_{ M^{ 2n}} \langle\sqrt{-1}(\bar \p_A \Lambda \p_A -\sqrt{-1} \p_A ^*\p_A) \phi , \psi\rangle - \int _{ M^{2n}} \langle \bar \p_A ^*\bar \p_A \phi, \psi \rangle $$
 $$ = \int_{ M ^{2n}} \langle\sqrt{-1}\Lambda \p_A \phi, \bar \p_A ^* \psi \rangle + \int _{ M^{2n}} \langle \p_A \phi, \p_A \psi\rangle
 - \int_{ M^{2n}} \langle\bar \p_A \phi, \bar \p_A \psi\rangle.\leqno (3.14)   $$
 
 Using  (3.13.1) we get Corollary 3.13.2 immediately from (3.14).
  
 \medskip
  
3) 
Using the Hodge-K\"ahler identities (3.1) and (3.2), we get
$$ -\sqrt{-1} \triangle ^{\p} _A  = \p_A ( \Lambda \bar \p_A - \bar \p_A\Lambda) + ( \Lambda \bar \p_A -\bar \p_A \Lambda) \p _A $$
$$ = \p_A \Lambda \bar \p_ A - \p _A \bar \p _A\Lambda  + \Lambda \bar \p_A \p _A - \bar \p_A \Lambda \p_A\leqno (3.15)$$
In the same way we have
$$\sqrt{-1} \triangle _A^{ \bar \p} = \bar \p _A ( \Lambda \p _A - \p_A \Lambda) + (\Lambda \p_A -
\p _A \Lambda) \bar \p _A $$
$$ =\bar \p _A \Lambda\p_A - \bar \p_A \p_A \Lambda + \Lambda\p_A \bar \p_A - \p_A \Lambda \bar \p_A . \leqno (3.16)$$
Using the identities
 $$ -(\p_A \bar \p_A + \bar  \p_A  \p_A) =  -F^{1,1} _A\wedge$$ 
 we get from (3.15) and (3.16)
$$ -\sqrt{-1}( \triangle _A^{\p}  - \triangle_ A ^{\bar \p}) = - [F_A^{1,1}\wedge, \Lambda].$$
 which yields (3.13.3) immediately.
 
 \medskip
 
4) We have
$$ \triangle _A ^d =  (\p_A + \bar  \p_A ) (\p_A^* + \bar \p_A ^*) + (\p_A^* + \bar \p_A ^*)
 (\p_A + \bar \p_A)$$
 $$ = \triangle  ^ \p_A +\triangle ^ {\bar \p} _A   + (\p_A \bar \p_A ^*  + \bar \p_A \p_A ^* + \p _A ^* \bar \p_A + \bar \p _A^* \p_A)\leqno (3.17)  $$
 Using the Hodge-K\"ahler identity (3.1), and   replacing $\p_A \p_A$ by $ F^{2, 0}_A \wedge$, we get
 $$( \p _A \bar \p_A ^* + \bar \p_A^* \p_A) = -\sqrt{-1} \p_A (\Lambda\p_A - \p_A \Lambda) -\sqrt{-1}( \Lambda \p_A - \p_A \Lambda)\p_A =\sqrt{-1} [  F^{2,0}_A \wedge, \Lambda].\leqno (3.18)$$
 Similarly
 $$(\bar \p _A  \p_A ^* +  \p_A^*\bar \p_A) = -\sqrt{-1} \bar \p _A ( \bar \p _A \Lambda - \Lambda \bar \p _A)  -\sqrt{-1}  ( \bar \p _A \Lambda - \Lambda \bar \p _A) \bar \p _A= -\sqrt{-1}  [ F^{0,2}_A \wedge, \Lambda].\leqno (3.19)$$
 Using corollary 3.13.3,  we get from (3.17), (3.18), (3.19)
 $$ \triangle  ^d _A = 2 \triangle ^{\bar \p} _A  - \sqrt{-1} [ F ^{ 1,1}_A\wedge, \Lambda] + \sqrt{-1} [ F^{0,2}_A\wedge, \Lambda] - \sqrt{-1} [ F^{0,2} _A\wedge, \Lambda]$$
 which yields  (3.13.4) immediately.
 
\QED

 \medskip
 
 {\bf 3.20. Remark.} Clearly  (3.13.2) follows directly from (3.13.3).  Furthermore, taking into account (2.2), we
 conclude that   all the formulas
 in  Corollaries 3.15  are valid, if we replace   bundle $E$ by bundle $End_JE$.

 \medskip
 
 Using Corollary 3.13.1, we observe that  a connection $A$ over a compact K\"ahler manifold is Yang-Mills bar, iff $\Lambda \p_A F^{0,2}  _A = 0$.
 We shall call a connection $A$ almost holomorphic, if $\p_A F^{0,2} _A  = 0$. Using the  Bianchi identity $\bar \p_A  F^{0,2} _A = 0 $, we  get that
 $\p _A F ^{ 0,2} = 0$, iff  $ d _A  F ^{ 0,2}_A = 0$.   Since $F^{2,0} _A = - ( F ^{0,2} _A) ^*$,  we observe that  $d_A F ^{0,2} _A = 0 $, iff $d_A F^{2,0} = 0$.
 Using the Bianchi identity
 $d_A F_A = 0$, we  observe that  $A$ is almost holomorphic, iff $d F^{ 1,1}_A = 0$. If $F^{1,1} _A = 0$ we shall call $A$ almost flat  holomorphic connection.  
 
  If dimension of $M$ equals 4,  it is easy to check that
 $$ \Lambda \p _A F ^{0,2} _A =  0  \LLR \p_A F^{0,2} _A = 0.$$
Thus  any Yang-Mills bar connection  over $M^ 4$  is an almost holomorphic connection.
 
 \medskip
 
 {\bf 3.21. Existence of almost holomorphic connections.}  Let  $T ^4$ be a   2-dimensional complex
 torus with  coordinates $z_1 = x_1 + \sqrt{-1} y_1$, $ z_2 = x_2 +\sqrt{-1} y_2$.
 Let $L$ be a complex line bundle whose Chern class is  represented by
 the cohomology class $c_1$ of $dz_1 \wedge dz_2 + d\bar z_1 \wedge d\bar z_2$. 
 Let $A$ be a unitary connection  of $L$. Then $F_A =  \sqrt {-1} (dz_1 \wedge dz_2 + d\bar z_1 \wedge d\bar z_2) + \sqrt{-1}d \alpha$, where $\alpha \in \Om^1 (T^4)$. The new connection $ A' = A -\alpha$ has  the curvature $\sqrt{-1} (dz_1 \wedge dz_2 + d\bar z_1\wedge  d\bar z_2)$, whose    component $F^{1,1}_{ A '} $  vanishes.  Thus $A'$ is an almost flat  holomorphic connection. We observe that by the  Hodge theorem $L$ carries
 no holomorphic structure. 
 
 The same argument  provides us  a differential-geometric proof of the Hodge conjecture for  Hodge classes
 of  dimension 2.
 
 \medskip
 To get  an almost holomorphic connection  in vector bundles
 of higher dimension we can take the sum of line bundles or a tensor product of a complex  line bundle
 with a holomorphic vector bundles.
 
  \medskip
  
  In the  next section we shall show that  if $M^{2n}$ is a K\"ahler manifold of positive 
  Ricci curvature, then any almost holomorphic connection is a holomorphic connection (Theorem 4.25), in particular
  any almost flat holomorphic connection is a flat connection.
  
  \medskip
  
  In general,   the Hodge theory implies that  on any  Hermitian complex line bundle over
  a K\"ahler manifold there is a  Yang-Mills bar connection which realizes the infimum  of the Yang-Mills bar functional.

  \medskip

 \section{Yang-Mills bar equation over  compact K\"ahler manifolds  of positive  Ricci curvature}

Suppose that $A$ is a  unitary connection on a Hermitian vector bundle $E$ over
a K\"ahler manifold $M^{2n}$. 
Let $D$ be the Levi-Civita connection on $T  ^*M^{2n}$:
$$D : \Om ^1 (M^{2n}) \to \Om ^1(M^{2n}) \otimes T ^* M^{2n}.$$
The connection $D$ extends  $\C$-linearly  to a  connection  also denoted by $D : \Om ^1 _\C (M^{2n}) \to \Om ^1_\C (M^{2n}) \otimes_\C  T_\C ^*M^{2n} = _\R  \Om ^1 _\C (M^{2n}) \otimes _\R  T ^* M^{2n}$.
Since $M^{2n}$ is K\"ahler,  we have
$D_v ( \phi \pm \sqrt{-1} J \phi )  = D_v (\phi) \pm \sqrt{-1} J D_v (\phi)$ for all $v \in T^*_\C M^{2n}$ and for all $\phi \in \Om^{0,1}  (M^{2n} )$.
It follows that  $D (\Om ^{0,1} (M^{2n})) \subset \Om ^ {0,1} (M^{2n}) \otimes_\C  T_\C ^* M^{2n}$, and iterating we have  $D (\Om ^{0,p} (M^{2n})) \subset \Om ^ {0,p} (M^{2n}) \otimes_\C  T_\C ^* M^{2n}$ for all $p$.
Now we denote by $\bar D$ the composition
$\pi ^{ 0,1} \circ D : \Om ^ {0,p} (M^{2n})  \to \Om ^{0,p} (M^{2n}) \otimes_{\C}  T ^{0,1} M^{2n}$, where
$\pi^{0,1}$ is the projection  to the  corresponding component with $(0,1)$-forms. Clearly for
all $\phi \in \Om ^{0,p}( M)$ the following formula holds
$$D_{v ^{0,1} } ( \phi) = \bar D_{ v^{0,1}}  (\phi),\leqno (4.1)$$
where $v ^{0,1}$ denotes the $(0,1)$-component of $v$: $ v ^{0,1} = (1/2) ( v + \sqrt{-1} J v)$.
Similarly, we shall use later the notation $v ^{1,0}  = (1/2) ( v - \sqrt{-1} J v)$.

 Combining  $\bar D$ with $\bar \p_A : \Om (E) \to \Om ^{0,1} (E)$, we define  the  following  partial connection
   $$\bar \nabla _A : \Om ^ {0,p} (E) \to  \Gamma ( E\otimes_\C\Lambda ^{0,p} T ^*_\C M^{2n} \otimes_\C  T ^{ 0,1} M^{2n}).$$
 In view of (4.1) we have
 $$\bar \nabla _A  = \pi ^{0,1} \circ \nabla A_{|\Om ^{0,p} (E)},$$
 where  $\nabla _A$ is the tensor product of $ d_A$ and $D$, which preserves the  natural induced metric on the
 bundle $E\otimes_\C \Lambda ^ p  T ^*_\C M^{2n} $
:
 $$\nabla _A : \Om ^p_\C (E) \to \Gamma (E\otimes_\C \Lambda ^ p  T ^*_\C M^{2n} \otimes_\C  T ^* _\C M^{2n}).$$
  In view of (4.1) we  also have
 $\nabla _A (\Om ^{0,p} (E) ) \subset \Om ^{0,p} (E) \otimes_\C  T^* _\C M^{2n}$.
 
 We shall  use the following notation. For any element $\phi \in \Om ^p ( E)$
 the expression $\phi _{  v_1,\cdots, v_p}$ denotes the value of $\phi$  at $(v_1, \cdots , v_p) \in \Lambda ^p (T_* M^{2n})$.
 
 Now we define a basic zero  order  operator $\Rr ^{ A} : \Om ^1_\C ( End_J E)\to \Om ^1_\C  (End_J E)$ by setting
 $$\Rr ^A (\phi)_X = \sum _{ j =1} ^{2n}[ (F _A) _{ e_j, X }, \phi _{e_j} ]\in End_J E,\leqno(4.2)$$
 where $(e_1, \cdots, e_{ n +k } = J e_k, \cdots,  e_{2n} )$ is a unitary basis  of the tangent space $T_x M^{2n}$ at the point
 $x$ in question. We also consider (as before) $F_A$ as an  element in $ \Om ^ 2_\C ( End _J ( E))$.

 Recall that the Ricci transformation $Ric : T_x M^{2n} \to T_x M^{2n}$ is defined by
  $$ Ric (X) = \sum _{ j =1} ^{2n} R_{ X, e_j} e_j,$$
 where $R$ denotes the curvature tensor of the Levi-Civita connection  on the tangent
 space $TM^{2n}$.  We denote by  the same $Ric$ the $\C$-linear extension of
$Ric$   from $T_xM^{2n}$ to  $(T_xM ^{2n}) _\C$. 

 We modify  this transformation by setting
 $$ Ric ^- (X): =  \sum _{ j =1} ^{2n} R_{ X, e_j} e_j ^{0,1}\in 
( TM ^{2n}) _\C .$$
Since $ J \circ R  = R \circ J$, we have  $ Ric ^- (X) = \pi ^{0,1}\circ Ric (X)$. Here $\pi ^{0,1}$ denotes
the projection on the $(0,1)$-component.
 
 Given $\phi \in \Om^1 (End_J E)$ we define    a new 1-form $\phi \circ Ric \in \Om ^1 (End _J E)$ by requiring  that
 for all $X \in  TM$ we have
 $$ (\phi \circ Ric )_X : = \phi _{ Ric (X)}.$$
 We also define  $\phi \circ Ric  ^- \in \Om  ^1 (End _J E)$ by requiring that  for all $X \in T M$ we have
 $$  (\phi\circ Ric^-) _ X: = \phi _{ Ric^- (X)}.$$
 
 If $\phi \in \Om ^{0,1} (End_JE)$, it is easy to see that $\phi \circ Ric ^- = \phi \circ  Ric$.

 \medskip
 
 {\bf 4.3. Lemma.} {\it   Suppose that  $(E, h) $ is a  Hermitian vector bundle provided with a unitary connection $A$. We have the following simple formulas for any $\phi \in \Om ^{0,p} (E)$ and  for arbitrary  $(0,1)$-vectors $X_i$
 $$(\bar \p_A \phi)_{ X_0, \cdots, X_p}  = \sum _{ k =0} ^ p (-1) ^ k( (\nabla _A)_{X_k}\phi ) _{ X_0, \cdots, \hat  X_k, \cdots X_p},\leqno (4.3.1)$$
 $$(\bar \p _A  ^ *  \phi) _{ X_1, \cdots , X_{p-1}} =  -\sum_{j =1} ^ {2n }( (\nabla _A)_{e_j ^ {1,0}} \phi ) _{e_j ^{0,1}, X_1, \cdots , X_{p-1}},\leqno (4.3.2)$$
where $ (e_1, \cdots , e_{n +k } = J e_{ k}, \cdots, e_{ 2n})$ is  an unitary frame  at a given point.}

\medskip


 {\it Proof.} First we extend  a well-known formula for  real forms (see e.g. \cite{B-L1981}, (2.12), (2.13)) to complex forms $\phi\in   \Om ^k_\C ( E)$ and
 $ X_i \in T^*_\C M ^{2n}$:
 $$( d_A \phi)_{ X_0, \cdots, X_p}  = \sum _{ k =0} ^ p (-1) ^ k(( \nabla _A)_{X^k}\phi )  )_{ X_0, \cdots, \hat  X_k, \cdots X_p},\leqno (4.4)$$
 
Formula (4.4)  holds, since  it holds for all real  forms $\phi\in \Om ^k(E) \subset \Om ^{k} _\C (E)$ and for all
$X_i \in T_* M ^{2n}$, and because both  LHS and RHS of (4.4) 
are $\C$-linear w.r.t.  to variables $\phi$ and $X_k$.

By definition  the LHS of (4.3.1) equals the LHS of (4.4) and clearly the RHS of (4.3.1) equals the RHS
of (4.4). Hence we get (4.3.1).

Now let us prove (4.3.2). For $ \phi \in  \Om ^{0,p} (E)$ and  for  a set of a $(1,0)$-vector $ X_0$  and $(0,1)$-vectors  $X_i$,  $1\le i \le p$, using (4.4), we have
$$  (\p _A  \phi)_{X_0 , X_1\cdots, X_{p}} = \sum _{ k =0} ^ p (-1) ^ k( (\nabla _A)_{X_k}\phi ) _{ X_0, \cdots, \hat  X_k, \cdots X_p},\leqno (4.5)$$
since LHS of (4.5) coincides with  the value $(d_A \phi) _{X_0 , X_1\cdots, X_p}$.  Since  $\phi\in \Om ^{0,p}(E)$, we get
$$\sum _{ k =0} ^ p (-1) ^ k( (\nabla _A)_{X_k}\phi ) _{ X_0, \cdots, \hat  X_k, \cdots , X_p} = ((\nabla _A)_{X _0} \phi) ( X_1, \cdots, X_p).\leqno (4.6)$$
Thus we get
$$ (\p _A \phi) = \sum_{ i = 1} ^ n dz_i \wedge  (\nabla _A )_{ e _i ^{ 1,0}  } \phi .\leqno (4.7)$$
Now  using the  K\"ahler identity  $\bar \p _A ^ * = - \sqrt{-1} \Lambda  \p_A $, we get from (4.7)
$$ (\bar \p_A  ^ * \phi) _{X_1, \cdots , X_{p-1}}  ={-1\over 2} \sum_{ k=1} ^n\sum_{j = 1} ^ n [ \bar i_k i_k dz_j\wedge((\nabla_A) _{e_j^{1,0}}\phi)]_{X_1, \cdots, X_{p-1}}$$
$$ = -\sum_{ j =1} ^ n[ \bar i_j ((\nabla _A)_{ e_j ^{ 1,0} }  \phi)] _{ X_1, \cdots ,  X _p}.\leqno (4.8) $$
 Clearly, the  last term of (4.8) equals the RHS of (4.3.2).  This completes the  proof of Lemma 4.3.
\QED
\medskip

The following Proposition is a complex analogue of Theorem 3.2 in \cite{B-L1981}.
\medskip

{\bf 4.9.  Proposition.} {\it For   any $\phi \in \Om ^{0,1} (End_J E)$  the following identity holds}
$$\triangle ^{\bar \p}_A \phi = \bar \nabla _A^* \bar \nabla  _A  (\phi) + \phi\circ Ric  + \Rr ^{ A} (\phi).\leqno (4.9.1)$$

\medskip

{\it Proof.}   Let $ X \in T_x ^ {0,1} (M^{2n})$.
We extend  $X$ locally on $M^{2n}$ so that $ DX  (x) = 0$.
We also extend the unitary frame $\{ e_1, \cdots, e_{ n+k} : = J e_k, \cdots, e_{2n}\}$ locally so that $  De_i (x) = 0$. Using (4.3.1) and (4.3.2), and taking into  account $(J e_j) ^{ 0,1}  = -\sqrt {-1}  e_ j ^{ 0,1} $,   $(J  e_j) ^{ 1,0} =   \sqrt {-1} e_j ^{1,0}$, we get at the point $x$
$$(\bar \p _A  \bar \p _A ^  * \phi)_X = (\nabla _A)_X \{\bar \p _A^* \phi \}   = - (\nabla _A)_{X} \{ \sum _{ j =1} ^{2n} [  (\nabla _A)_{e_j} \phi] _{e_j ^{0,1} }\} $$
$$ = - \sum _{ j =1}^{2n} [( \nabla _A)_X   (\nabla _A) _{e_j}\phi ]_{ e _j ^{0,1} }.\leqno(4.10)$$

\medskip

$$ ( \bar \p_A  ^ * \bar \p _A \phi) _X = -\sum_{ j=1 } ^{2n} \{(\nabla _A)_{e_j} (\bar \p _A \phi)\} _{ e_j ^{0,1} , X}$$
$$ =-\sum_{ j=1} ^{2n}(\nabla_A)_{e_j} \{ [(\nabla_A)_{e_j^{0,1}} \phi  ]  _X - [(\nabla_A)_X \phi] _{ e ^{ 0,1}_j} \}$$
$$ = -\sum_{j=1} ^{2n} \{ [(\nabla _A)_{e_j} (\nabla _A)_{e_j ^{0,1}} \phi]_X  - [(\nabla _A)_{e_j} (\nabla _A)_X \phi]  _{ e _j ^{0,1} }\}. \leqno(4.11)$$

Summing (4.10) and (4.11), we get

$$(\triangle ^{\bar \p}_A  \phi)_X  = -\sum_{ j =1} ^{2n} \{ [(\nabla _A)_{e_j} (\nabla _A)_{e_j ^{0,1}} \phi]_X + \sum _{ j =1} ^{2n} ( R ^{A}_{ e_j, X} \phi) _{ e^{0,1}_j}\}.\leqno (4.12)$$

Here we denote by $ R ^A$  the curvature of the tensor product connection on the bundle
$T ^*_\C M ^{2n} \otimes_\C  End_J E = ( T^*M ^{2n}\otimes_\R End_J E)_\C$. This curvature coincides with the one  on $T^*M^{2n} \otimes _\R End _J E $, if
we consider $\Om ^2 (T^*M \otimes_\R End _J E)$ as a subspace in $\Om ^2_\C (T^*_\C M^{2n} \otimes_\C End_J E)$.
Now we observe that for $\psi \in \Om ^{0,1} (End_J E)$ we get
$$\langle-\sum _{ j =1} ^{ 2n} [(\nabla _A)_{e_j} (\nabla _A)_{e_j ^{0,1}} \phi] , \psi\rangle= $$
$$  -\sum_{j=1} ^{2n}[(\nabla_A)_{e_j} \langle(\nabla_A)_{e_j ^{0,1}} \phi, \psi\rangle - \langle (\nabla_A)_{e_j ^{0,1}} \phi , (\nabla_A)_{e_j} \psi\rangle ].\leqno (4.13)$$

We define a 1-form  $\sigma$, depending on $\phi$ and $\psi$,  on $M$  by
$$\sigma(X) := \langle(\nabla_A)_{X ^{0,1}} \phi, \psi\rangle .$$
Then  $$- \sum_{j=1} ^{2n}(\nabla_A)_{e_j} \langle(\nabla_A)_{e_j ^{0,1}} \phi, \psi\rangle(x)  = (- d  ^* \sigma) (x),\leqno(4.14)$$
and
$$  \langle (\nabla_A)_{e_j ^{0,1}} \phi , (\nabla_A)_{e_j} \psi \rangle  = \langle \bar \nabla _A \phi , \bar \nabla _A \psi\rangle
+ \sum _{ i = 1} ^{ 2n}\langle (\nabla _A)_{ e_j ^{ 0,1}} \phi, (\nabla _A) _{ e_j ^{ 1,0}} \psi\rangle.$$
Since $(J  e_j) ^{ 0,1}  = -\sqrt {-1}  e_ j ^{ 0,1} $ and $(J  e_j) ^{ 1,0} =   \sqrt {-1} e_j ^{1,0}$, we get
$$\langle(\nabla_A)_{  e _i ^{ 0,1} } \phi, (\nabla _A) _{ e_j ^{ 1,0}} \psi\rangle  +\langle(\nabla_A)_{ (J e _i) ^{ 0,1}}  \phi, (\nabla _A) _{ (Je_j )^{ 1,0}} \psi\rangle  = 0$$
$$ \LRA \sum _{ i = 1} ^{ 2n}\langle (\nabla _A)_{ e_j ^{ 0,1}} \phi, (\nabla _A) _{ e_j ^{ 1,0}} \psi\rangle= 0$$
 $$\LRA  \langle (\nabla_A)_{e_j ^{0,1}} \phi , (\nabla_A)_{e_j} \psi \rangle  = \langle \bar \nabla _A \phi , \bar \nabla _A \psi\rangle.\leqno (4.15)$$
 From (4.13), (4.14), (4.15) we get
 $$\int_{ M^{2n}} - \langle (\nabla _A)_{ e_j} ( \nabla _A)_{ e_j ^{0,1}} \phi, \psi \rangle = \int_{ M^{2n}}-d ^* \sigma +\int _{ M^{2n}} \langle \bar \nabla _A \phi, \bar\nabla _A \psi \rangle.\leqno (4.16)$$ 
 Next we  have
 $$ (R ^A _{ e_j, X} \phi) _{ e_j ^{0,1}} =  (F _A) _{ e_j, X } \phi _{ e_j ^{0,1}} -
  \phi ( R _{ e_j , X } e_j ^{0,1} ). \leqno (4.17)$$

  Clearly Proposition 4.9 follows from (4.12), (4.13) and (4.17).

\QED

\medskip
 
Denote by $  \Rr ^A $ the following  linear operator  $:\Om ^{0,2} (End_J E)\to \Om ^{0,2} (End_J E)$ such that  for all
$(0,1)$-vectors $X, Y$ we have
$$ (\Rr ^A  ( \phi))_{ X, Y }  = \sum_{ j =1} ^{2n}  \{  [ (F _A) _{ e_j, X} , \phi_{e_j, Y } ]- [ (F_A) _{ e_j , Y}, \phi _{ e_j , X } ]\}.\leqno(4.18)$$
 
 We also associate to  each $\phi \in \Om^{0,2} ( End_J E)$ a new
 $(0,2)$-form $\phi \circ ( Ric \wedge I) \in \Om^{0,2} ( End_ J E)$  by setting
 $$( \phi \circ ( Ric \wedge I)) _{ X, Y } : =  \phi  ( Ric (X), Y) - \phi ( Ric  (Y), X ), $$
 
\medskip

{ \bf 4.19. Proposition}.  For  any $\phi \in \Om ^{0,2} (End_J E)$  the following identity holds
$$\triangle  ^{\bar\p} _A \phi  = \bar \nabla _A ^*\bar  \nabla_A  \phi + \phi \circ ( Ric   \wedge I )  +  \Rr ^A ( \phi).\leqno (4.19.1)$$

\medskip

{ \it  Proof.} (Cf.  Theorem 3.10 in \cite{B-L1981}.) We use the notations  $X, Y, e_1, \cdots , e_n$ as in the proof of Proposition 4.9. Then at the  point $x$ and for  $(0,1)$-vectors $X$ and $Y$  we have
$$ (\bar \p _A \bar \p _A  ^*  \phi ) _{ X, Y} = ((\bar \p_A )_X \bar \p _A ^*  \phi) ) _Y
- ( (\bar \p_A ) _Y \bar \p_A ^* \phi)_X $$
$$ =-  (\nabla_A)_X \{  \sum _{ j= 1} ^{2n} ((\nabla _A)_{ e_j} \phi) _{ e_j ^{0,1} , Y }  \}
+ ( \nabla_A ) _Y  \{  \sum _{ j =1} ^{ 2n}  ((\nabla _A ) _{ e_j} \phi) _ { e_j ^{0,1}, X} \} $$
$$  = -\sum _{j = 1} ^{ 2n} \{[( \nabla _A)_X (\nabla _A)_{ e_j} \phi] _{ e_j ^{0,1}, Y} 
- [( \nabla _A)_Y (\nabla _A)_{ e_j} \phi] _{ e_j ^{0,1}, X} \}.\leqno  (4.20)$$

We also have
$$(\bar \p _A ^* \bar \p _A \phi)_{ X, Y} =  -\sum _{ j =1 } ^{ 2n} ( (\nabla _A )_{ e_j} \bar \p_A \phi)_{ e_j ^{ 0,1}, X, Y }  $$
$$ = -\sum _{ j = 1} ^{2n}( \nabla_A)_{ e_j} \{( (\nabla _A)_{ e_j^{0,1}}  \phi ) _{ X, Y } +
((\nabla _A)_Y \phi )_{  e_j^{ 0,1} , X } + ((\nabla _A)_X \phi )_{Y,   e_j^{ 0,1} }\} $$
$$ = -\sum_{ j = 1} ^{2n} \{[( \nabla_A)_{e_j} (\nabla _A)_{ e_j^{ 0,1} } \phi] _{ X, Y}
 + [( \nabla _A )_{ e_j}  (\nabla _A) _ Y \phi] _ { e_j ^{0,1} , X }  - [ ( \nabla _A)_{ e_j} 
 ( \nabla _A) _X \phi ] _{ e_j^{0,1} , Y} \} . \leqno (4.21)$$

Summing  (4.20) and (4.21) we get
$$( \triangle  ^{ \bar \p } _A  \phi )_{ X, Y } = -\sum_{ j = 1} ^{2n} [( \nabla_A)_{e_j} (\nabla _A)_{ e_j^{ 0,1} } \phi] _{ X, Y} + \sum _{ j = 1} ^{2n}  \{[ R^A _ { e_j,  X } \phi ] _{ e_j ^{0,1} , Y} -  [ R ^A _{ e_j, Y } \phi] _{ e_j^{ 0,1} , X} ] \}. \leqno (4.22)$$

As in the proof of Proposition 4.9 (see (4.13)) we have for $\psi \in \Om ^{0,2} (End_JE)$
$$\int_{ M^{2n}}\langle -\sum _{ j = 1} ^{2n}( \nabla_A)_{ e_j} (\nabla _A)_{ e_j^{0,1}}  \phi , \psi\rangle
= \int_{ M^{2n}} \langle \bar \nabla _A \phi, \bar \nabla _A \psi \rangle.\leqno (4.23)$$

We  use the following identity
$$  ( R  ^A_{ X, Y } \phi) _{ Z, W } =  [ (F_A)_{ X, Y }, \phi_{ Z, W } ] -
\phi ( R_{ X, Y } Z, W ) - \phi  ( Z, R _{ X, Y } W ) $$
and combining with (4.23)  to rewrite  (4.22) as follows
$$( \triangle  ^{ \bar \p } _A  \phi )_{ X, Y } = ( \bar \nabla _A ^* \bar  \nabla _A  \phi) _{ X, Y } + \Rr ^A(\phi) _{ X, Y }  +  \phi ( Ric^- (X), Y)  $$
$$ - \sum _{ j =1} ^{2n}\phi _{ e_j ^{ 0,1}, R _{ e_j, X } Y} - \phi  ( Ric ^- (Y), X ) + \sum_{ j =1} ^{ 2n}  \phi _{ e_j ^{0,1}, R_{ e_j, Y } X } .\leqno (4.24)$$
Using the Bianchi identity
$$-  R_{e_j ,  X } Y  - R_ { Y, e_j } X =  R_{ X, Y } e_j, $$
and taking into account that  the following quantity vanishes for all $\phi \in \Om^{0,2}  ( End_J E)$  and for all $X, Y \in T ^{0,1} M^{2n}$
$$ (\phi \circ R)  _{ X, Y } : = \sum_{ j =1 } ^{2n} \phi ( e_j,  R _{ X, Y } e_j ),$$
because  $(J  e_j) ^{ 0,1}  =- \sqrt {-1}  e_ j ^{ 0,1} $ and $(J  e_j) ^{ 1,0} =   \sqrt {-1} e_j ^{0,1}$,
we  get Proposition 4.19  from (4.24)  immediately.\QED

\medskip




{ \bf 4.25. Theorem.} {\it   Let $M$ be  a compact K\"ahler manifold with positive  Ricci curvature. If $A$ is an almost holomorphic connection,
then $A$ is holomorphic.}

\medskip

{\it Proof.}  First let us prove the following formula  for $\phi \in \Om ^{ 0,2}( End_J E )$. 
$$ \Rr ^ A( \phi) =  -  \sqrt{-1} \{\Lambda F^{1,1}_A\wedge \phi - (\Lambda F ^{1,1}_A ) \phi\}
:=\bar R (A) \phi.\leqno  (4.26)$$
Let us  rewrite the expression in (4.19)  as follows
 $$ \Rr  ^A (\phi)  = \sum_{ 1\le k  < l \le n} \sum _{ j = 1} ^{2n}\{ [(F_A)_{e_j,e_k ^{0,1}}, \phi _{ e_j, e_l ^{0,1}}] - [(F_A) _{ e_j , e_l ^{0,1}},\phi_{ e_j, e_k ^{0,1}}]\} d\bar z_k d \bar z_l .\leqno (4.27)$$
 
 We shall use the following abbreviation. For any $\phi \in \Om ^{k,p} (End_J E)$ denote by
 $$\phi _{ i_1\cdots i_k, \bar j_1\cdots \bar  j_p}: = \phi ( e_{i_1} ^{1,0}, \cdots, e_{i_k} ^{1,0}, e_{j_1} ^{0,1},
 \cdots , e_{j_p} ^{0,1}).$$
 Since $\phi \in \Om^{0,2}(End_J E)$, we get from (4.27)
 $$\Rr ^A (\phi) =  \sum_{ 1\le k < l \le n}  \sum_{ j = 1} ^{2n}\{  [(F_A)_{j \bar k},\phi _{ \bar j\bar l}]-
 [(F_A) _{ j \bar l}, \phi_{\bar j\bar k} ] +  [(F_A)_{\bar j \bar k},\phi _{ \bar j\bar l}]-
 [(F_A) _{\bar j \bar l}, \phi_{\bar j\bar k} ]\}d\bar z_kd \bar z_l  .\leqno (4.28) $$
 Since $(J ( e_j)) ^{ 0,1}  = -\sqrt {-1}  e_ j ^{ 0,1} $ and $(J ( e_j)) ^{ 1,0} =   \sqrt {-1} e_j ^{1,0}$, we get from (4.28)
 $$\Rr ^A (\phi) = 2 \sum_{ 1\le k < l \le n}  \sum_{ j = 1} ^{n}\{ [(F_A)_{j \bar k},\phi _{ \bar j\bar l}]-
 [(F_A) _{ j \bar l}, \phi_{\bar j\bar k} ]\}\, d\bar z_ k d\bar z_ l.\leqno (4.29) $$

Now  expanding the  following  expression  in local coordinates

$$\sqrt{-1} \Lambda F ^{1,1} \wedge \phi = {1\over 2}\sum _{p = 1} ^n \bar i _p i _p \{ \sum _{ i,j} \sum _{ k < l}[(F_A) _{i\bar j}, \phi _{\bar k \bar l}] \, dz_id\bar z_j d\bar z_k d\bar z_l\} $$
= $$ \sqrt{-1}( \Lambda  F ^{1,1} _ A ) \phi  - 2 \sum_{ 1\le i \le n}\sum _{1\le j, l\le n}  [(F_A) _{ i \bar j } , \phi _{ \bar i\bar l}] d\bar z_j d\bar z_l,$$ 
and comparing it with  the RHS of (4.29),   we get
(4.26) immediately.  
\medskip

Now let  $A$ be a Yang-Mills bar connection.  Applying   (4.19.1) to $F^{0,2} _A$  and using (4.26) we get

$$ 0  = \int_{M^{2n} } \langle \bar \nabla _A   F^{0,2} _A, \bar \nabla  _A F^{0,2} _A \rangle 
+\langle  F^{0,2}_A\circ ( Ric \wedge  I), F^{0,2} _A \rangle + \int _M \langle \bar R(A) F ^{0,2}_A, F ^{0,2}_A\rangle.\leqno (4.30)    $$
Since $A$ is a Yang-Mills bar connection, differentiating (2.7.1), 
 we get 
$$ \langle (\Lambda F ^{1,1}_A) F ^{0,2}_A, F ^{0,2 } _A \rangle = 0.\leqno (4.31)$$

Now let $A$ be an almost holomorphic connection.  Using  (4.26), (4.30), (4.31), (3.13.2) (see also corollary 3.20), we get  $F^{0,2} _A =0$  immediately. \QED
 


\medskip

{\bf 4.32. Remark.} Theorem 4.25 implies that any  Yang-Mills bar connection on a compact 4-dimensional K\"ahler manifold of positive Ricci curvature is holomorphic.   
\medskip

 \section {Short time existence of a Yang-Mills bar  gradient flow over  a compact K\"ahler manifold}

 {\bf 5.1. Affine integrability  condition.}   The following identity holds for  any $\theta \in \Om(End_J E )$ and  any unitary connection $A$
 $$ \int _{ M^{2n}} \langle [ \theta, F ^{0,2}_A], F^{0,2} _A \rangle = -  \int_{ M^{2n} }\langle[ F^{0,2} _A, \theta], F^{ 0,2} _A \rangle. \leqno (5.2)$$
 We shall prove that at any point $x \in M^{2n}$ we have
 $$\langle[\theta, F^{ 0,2}_A], F^{0,2} _A \rangle =  - 2 \langle \theta , \Lambda \Lambda F^{0,2} _A \wedge F^{2,0} _A \rangle.\leqno (5.3) $$
 We write   $ \theta = \theta^ +  + \sqrt {-1} \theta ^{-}$ where $\theta^ +, \theta ^- \in  u _E$. In the same way  at a  fixed point $x \in  M^{2n}$ we can take  coordinates such that the K\"ahler metric $g$ has the form
 $g(x) =\sum dz_i \otimes d\bar z_i$.  We shall  write
 $$ F ^{0,2} _A = \sum_{1\le i< j\le n}(  F_{ ij} ^+  + \sqrt{-1} F^-_ { ij} ) d\bar z_i d \bar z_j,$$
 where $F ^{ \pm}_{ij} \in u _E$. Then $ F^{ 2,0} _A =  \sum _{ ij}(F^ +_{ij} - \sqrt{-1} F ^{-} _{ij}) dz_i dz_j$. Recall
 that $|| d\bar z_i d\bar z_j || ^ 2  = 4$.  A direct computation at a point $x$ shows
 $$\langle[ \theta, F ^{ 0,2} _A], F^{ 0,2} _A \rangle =\sum _{1\le  i< j\le n}  \langle [\theta ^{ -}, F^{+} _{ ij }] d\bar z_i d\bar z_j,  F ^{ -}_{ ij}  d\bar z_i d\bar z_j\rangle+$$
$$ \sum _{1\le i<j\le n}\langle- [\theta ^{-}, F ^{ -} _{ij}]d\bar z_i d\bar z_j , F ^{ +} _{ ij} d\bar z_i d\bar z_j\rangle = 8 \langle \theta^-, \sum _{1\le  i< j\le n}[ F^{+} _{ij}, F ^{ -} _{ ij}]\rangle. \leqno (5.4)$$
 Now we compute
$$\langle \theta, \Lambda \Lambda F ^{ 0,2}_A \wedge F ^{ 2,0} _A \rangle = -2\sum_{1\le i<j\le n} \langle \theta ^ {-},\Lambda \Lambda [ F^+_{ ij }, F ^- _{ ij}]d z_i d z_j d\bar z_i d\bar z_j\rangle   $$
$$   = - 4\sqrt{-1} \sum _{1\le i< j\le n} \langle \theta ^-,  \Lambda [ F^+ _{ij} , F^- _{ ij} ] (d  z_j d \bar z_j+ dz_i d\bar z_i)\rangle  =
 -16 \langle \theta ^- , \sum_{1<i<j\le n} [ F ^+_{ ij}, F^- _{ ij}]\rangle.\leqno(5.5) $$
Clearly (5.3) follows from  (5.4) and (5.5).

 Now  substituting  $[F ^{ 0,2} _A , \theta] = \bar \p _A \bar \p _A  \theta$ in the RHS of  (5.2),  and taking into account  (5.3), we get
 $$ \int_{ M^{2n}}  \langle \theta, 2\Lambda \Lambda F^{0,2} _A \wedge F^{2,0} _A \rangle =  \int _{ M ^{2n}}\langle \theta, \bar \p _A ^ * \bar \p_A ^ *  F ^{ 0,2} _A \rangle.\leqno (5.6)$$
 Thus we get the following identity
 $$\bar\p_A ^ * \bar \p _A ^ * F^{0,2} _A - 2\Lambda \Lambda F^{0,2} _A \wedge F^{2,0} _A = 0.\leqno (5.7)$$
 Define the following operator $ P_A : \Om ^{0,1} (End_J E)  \times \Om ^{ 0,1} (End_J E) \to \Om  (End_J E)$
 $$ P_A ( a) \phi : = \bar \p_{ A +a}  ^* \phi -   2 \Lambda \Lambda F ^{0,2} _{ A + a}\wedge  F ^{ 2,0}_ { A + a} . \leqno (5.8)$$
 Clearly $P_A (a) \phi$  is a differential operator  of order 1 in $a$ and order 1 in $\phi$.
 Moreover  $P_A(a) \phi$ is an affine differential operator w.r.t.  $\phi$, i.e.
 $P_A (a) \phi  = L _A(a) \phi + C_A (a)$, where $L_A(a)\phi$ is a linear differential operator
 w.r.t. $\phi$.
 By (5.6) we have $ P _A(a) \bar \p _{A+a}  ^ *F ^{ 0,2} _{ A + a}= 0.$  Thus  we shall call $P_A(a)$ an  affine  integrability condition 
 for  the differential operator $\bar \p_{A+a} ^ * F ^{ 0,2}_{A+a}: \Om ^{ 0, 1}( End_J E) \to  \Om ^{0,1} ( End_J E ).$
 
 \medskip
 
 {\bf  5.9. Proposition.} {\it   Let $ \xi  \in T_x ^ * M ^{ 2n} \setminus \{ 0\}$. All the eigenvalues of the eigenspace of
 the symbol  $\sigma_\xi D (-1) \p _{A+a} ^ *  F ^{0,2} _{ A + a} :\Om ^{ 0,1} (End_J E) \to \Om ^ { 0,1} ( End_JE) $  in Null $ \sigma _\xi P_A (a)$ are positive.   Hence the evolution equation
 $$ { d a \over dt} = - \bar \p _{ A + a}  ^ * F ^{ 0,2} _ { A + a} , \leqno (5.9.1) $$
 has  a unique smooth solution  for a short time which may depend on $a$.}
 
 \medskip
 
 {\it Proof.} Since $F ^{ 0,2} _{ A + a + th } = F^{ 0,2} _A + t \bar \p_{ A + a }\wedge  h +  t ^ 2 h \wedge h$ for
  $h \in \Om ^ {0,1} (End_JE)$, we  have the following expression for the   linearization of $\bar \p_{A+a} ^ *  F ^{ 0,2} _{A+a} $ at  point $a
  \in \Om ^{0,1} (End _J E)$ 
  $$ D_a (\bar \p_{ A + a } ^ *  F ^ {0,2} _ { A + a} ) (h) = \bar \p_{A+a} ^ * \bar \p_{ A + a}  h + \{ \text { terms of lower order } \}.\leqno (5.10)$$
   We may assume that $\xi  = d x_1$. 
  Then   a direct computation  using  the Hodge-K\"ahler identity $\bar \p _{ A + a }^* = - \sqrt{-1} \Lambda \p_{ A +a}$           and (5.10) shows
 $$-\sigma _\xi D_a ( \bar \p _ {A + a}^ *  F ^{ 0,2} _{A+a}) (  \alpha _1 d\bar z _1 , \cdots , \alpha _n d\bar z _n ) =  ( 0, \alpha _2  d \bar  z_2, \cdots,  \alpha _n d\bar z _n).\leqno (5.11)$$ 
 Clearly the linearization $D_{\phi} P_A (a)\phi$ with respect to the variable  $\phi$ is
 $$[D_{\phi} P_A (a) \phi ]h = {d \over dt}_{ |t = 0}\bar \p_{ A +a}  ^* (\phi + t h) -   2 \Lambda \Lambda F ^{0,2} _{ A + a}\wedge  F ^{ 2,0}_ { A + a} = \bar\p_{ A + a}^*  ( h).$$
   We note that this linearization does not depend on $\phi$.  A short computation shows
 $$\sigma _\xi  D _\phi P _A (a) ( \alpha _1 d\bar z_1, \cdots , \alpha _n d\bar z _n ) = \sqrt{-1} \alpha _1.\leqno (5.12)$$
 Now  (5.11) and (5.12)  imply the first statement of Proposition 5.9. The second  statement  follows  from  Hamilton's theory for
  evolution equation with integrability condition \cite{Hamilton1982}, Theorem 5.1, actually  from its slightly extended version in Theorem 6.6 below. 
 
 {\bf 5.13. Remarks}. 1.  By taking  derivative  of (2.7.1) in  the time $t$ we also get  (5.2) and hence  (5.7).  In the same
 way we can get  (5.2) (and hence (5.7))  as an infinitesimal consequence
 of the  non-canonical action  of the complex gauge group on the space of unitary connections w.r.t. a fixed  Hermitian metric on the bundle.
 
 2. It is likely that  $\bar\p_A ^* F ^{0,2} _A$
 also satisfies  an affine integrability  condition analogous to (5.8), if the ground manifold  $M^{2n}$ is  Hermitian but not
 necessary K\"ahler.
 
 \medskip

 \section{ Evolution equations with  affine integrability condition}
 
  In his work  \cite{Hamilton1982} Hamilton introduced the notion of  an evolution equation with integrability condition. Let us  rapidly  recall
 the  Hamilton concept from  section 5 of that paper.  We  try to  keep most of  Hamilton's notations in that paper,  which
 may have  quite different meanings  from that ones we used  in previous sections.
 
 \medskip
 
 We shall consider   an evolution equation
 $$ {df\over dt} = E (t), $$
 where $E(f)$ is a non-linear  differential operator  of degree 2 in $f$.  We suppose   that the values of  $f$ belong to an
 open set  $U$ in a vector bundle $F$ over a compact manifold $X$, and $E (f)$  takes its values in $F$  also.
 (For the case we are dealing in this note,  we shall take $U = F$. We shall write later, following Hamilton,  $f\in U$,
 meaning  that  the  values of $f$ belong to $U$.) Then $E$ is a smooth map
 $$E :  C^ \infty (X, U)  \subset  C ^\infty  (X, F)  \to C ^\infty (X, F)$$
 of an open  set in a Fr\'echet space to itself.
 
 We shall consider problems where some of the eigenvalues of  the symbol $\sigma  DE (f) \xi$ are zero. This happens when $E(f)$ satisfies  an integrability
 condition.
 \medskip
 
{\bf  6.1. Definition.} \cite{Hamilton1982} Let $g =  L (f) h: C^\infty (X, U) \times  C ^\infty (F) \to C ^ \infty (G)$ 
 be a differential operator  of degree  1  on sections $f \in U \subset  F$,  $h\in F$,  and $G$  another  vector bundle over $X$. We call $L(f) h$  the integrability condition
 for $E(f)$, if  the operator  $ Q(f)  = L(f) E(f)$ only has  degree at most one in $f$.
 
 \medskip
 
  Suppose that $L(f)h$ is an integrability condition  for $E(f)$. Taking a variation in $\tilde f$ we see that
 $$ L(f) DE (f) \tilde f  + DL (f) \{ E (f), \tilde f \}  =  D Q ( f) \tilde  f.\leqno (6.2) $$
 Since  $DQ(f) \tilde f$ as well as $L  (f) DE(f)  \tilde f$  only have degree 1 in $f$ the operator  $L (f) DE (f) \tilde f$ also have degree 1. 
 hence $\sigma L(f) (\xi)\sigma  DE (f) (\xi) = 0$. Therefore we get  
 $$ Im \, \sigma DE (f) ( \xi) \subset  Null \, \sigma  L (f) ( \xi). \leqno  (6.3)$$

 \medskip
 
 { \bf  6.4. Theorem} (\cite{Hamilton1982}, Theorem 5.1). {\it  Let $ df/ dt = E (f)$ be an evolution
 equation with integrability condition $L(f)$. Suppose that all the eigenvalues of the
 eigenspaces of $\sigma DE (f) (\xi) $ in $Null\,  \sigma L (f) (\xi)$ is positive. Then the
 initial value problem $f = f_0$
 at $t = 0$ has a unique smooth solution for a short
 time $ 0 \le t \le \eps $ where $\eps$ may depend on $ f_0$.}
 
 \medskip
 
 {\bf 6.5. Remark.}  Hamilton's notation in (6.2) indicates that $L (f) h$ is a linear w.r.t.
 $h$. (In fact, in section 4 of that paper Hamilton  stressed that $L(f) h$ is linear w.r.t. $h$.)  A  closer look at Hamilton's proof (see also  our proof of Theorem 6.6 below) shows that,
 the linearity of $L(f) h$ w.r.t. $h$ is important. We  shall call such integrability
 condition $L (f) h$  linear in the argument (and $f$  shall be considered as parameter).
 Now we shall call an  integrability condition $L (f) h $ an affine  integrability condition,
 if $ L(f) h = L_0 (f) h + A( f)$, where $L_0 (f) h$ is linear w.r.t. $h$. The   linearization
$(D_\phi L (f) h )  \tilde h = L_0 ( f) \tilde h$  does not depend on $h$.

 \medskip
 
 { \bf  6.6. Theorem.}   {\it  Let $ df/ dt = E (f)$ be an evolution
 equation with affine integrability condition $L(f)$:
 $L(f)h =L_0 (f) h + A ( f)$. Suppose that all the eigenvalues of the
 eigenspaces of $\sigma DE (f) (\xi) $ in $Null\,  \sigma L_0 (f) (\xi)$ is positive. Then the
 initial value problem $f = f_0$
 at $t = 0$ has a unique smooth solution for a short
 time $ 0 \le t \le \eps $ where $\eps$ may depend on $ f_0$.}
 
 \medskip

 {\it Proof of Theorem 6.6.} We follow  Hamilton's argument, replacing  $L(f) h$  in his proof
 by $L_0(f) h$ in some places, and  re-arranging  parameters which do not depend on $h$. To keep  our notations  as close as possible with
 those of Hamilton, we   denote by $DL$ the derivative of $L (f) h$ w.r.t the parameter $f$. We divide  the proof in  3 steps.
 
 \medskip

STEP 1. {\it Reduction of Theorem 6.6  to a version of the Nash-Moser inverse  function theorem.}

 In this step we reduce   Theorem 6.6 to the following 
 
\medskip

{\bf  6.7. Lemma}. { \it  Suppose that  $\bar f$ is a solution of   the perturbed  evolution equation by a term $\bar h(t,x)$
$$ {d\bar f(t,x)\over dt}  =  E ( \bar f(t,x))  + \bar h(t,x), $$
$$ \bar f (0,x) = \bar f_0(x)$$
over the interval $0 \le t \le 1$.
Then  for any  $f_0$ near $\bar f _0$ and $ h$ near  $\bar  h$ there exists a unique  solution of the    perturbed equation
$$ {df(t,x)\over  dt } = E (f(t,x)) + h(t,x) , $$
$$ f(0,x)  =  f_0(x)$$
over the interval  $0 \le t \le 1$.}

\medskip

Now we explain how   to get Theorem  6.6 from  Lemma 6.7. Let $\bar f(t, x)$   be any   function satisfying
$$ {d \bar f ( t, x) \over dt }_ { | t = 0}   =   E (  f ( 0, x)), $$
$$\bar f ( 0, x) =  f_0 (x).$$

Set 
$$ \bar h(t, x) :=  {d\bar  f(t, x) \over dt } - E ( \bar f ( t, x)).$$
Then $ \bar h (0,  x) = 0.$

Since $X$ is compact,  for any $ \delta > 0$ there exist a number  $ \eps > 0$  and  a function $ h (t, x)$ such that $H(t,x)$ is  $\delta$-close to  $\bar h (t, x)$ and
moreover $h (t, x ) = 0$ for  a short time $ 0 \le t \le \eps$.  Applying Lemma 6.7  to the pair $( \bar h, h)$ we conclude that
the equation
$${df(t,x) \over dt } =  E (f(t,x)) + h(t,x) , $$
$$f( 0, x)  = f _0 (x) $$
has  solution up to time $\eps$. This solution in the interval $( 0, \eps)$ is a solution  of our  original equation  in that time interval.
This completes the first step.

\medskip

STEP 2. { \it Reduction  of Lemma 6.7 to a case  of  a  weakly parabolic linear system  of (6.14.1) and (6.14.2).}
We can  apply the Nash-Moser  inverse function theorem to the operator
$$ \Ee : C^\infty ( X \times [0,1] , F ) \to  C^\infty  ( X \times [0,1], F) \times C^\infty ( X, F), $$
$$ \Ee (f) = ( df/ dt - E (f), \, f | \{ t = 0 \}).$$

Its derivative is the operator 
$$ D \Ee ( f) \tilde f = ({ d\tilde f \over  dt} - DE (f) \tilde  f, \tilde f| \{ t = 0 \}).$$
We must  show that the linearized  equation  
$$ d\tilde f / dt - DE (f) \tilde f = \tilde h \leqno (6.8)$$
has a unique solution for the initial value problem $\tilde f = \tilde f_0$ at $t = 0$, 
and verify that the solution  $\tilde f$ is a smooth tame function of $\tilde h$ and $\tilde f_0$.

 We make  the substitution  $ \tilde g = L (f) \tilde f$.  Then $\tilde g $ satisfies the evolution equation
 $${ d \tilde g \over dt } = L_0( f ) { d \tilde f \over d t} + DL (f) \{ \tilde f,  { d f\over d t}\} .\leqno (6.9)$$
 
 Now differentiating the integrability condition $L (f) E(f)  = Q (f)$ we get
 $$L_0(f) DE (f) \tilde f  = -DL (f)  \{ E (f), \tilde f \}  +  DQ (f) \tilde f. \leqno (6.10)$$
 
 Substituting $d\tilde f / dt   =  DE (f) \tilde f + \tilde h$ from (6.8) into (6.9) and taking into account  (6.10) we rewrite (6.9) as follows
 $${d\tilde g \over dt } - M (f) \tilde f = \tilde k, \leqno (6.11)$$
 where $\tilde k = L_0 (f) \tilde h$ and 
 $$M (f) \tilde f =  DL (f) \{ \tilde f, { df \over dt  }   \} -  DL (f) \{ E (f), \tilde f \} + DQ (f) \tilde f =$$
 $$ \stackrel{(6.10) } {=}  DL ( f)  \{\tilde f,  { df \over  dt } \}  + L_0 (f) DE (f) \tilde f. \leqno (6.12)$$ 
 is a  linear  differential operator in $\tilde f$ of degree $1$  whose coefficients depend smoothly  on
 $f$ and its derivatives.
 
  If we choose a measure on $X$ and inner product on the  vector bundle $F$ and $G$, we can form  a  differential operator  $L_0 ^* (f) g  = h$ of degree 1 in
  $f$ and  $g$  which is the adjoint of $L_0(f)$.  Let us write
 $$ P (f) h : = DE (f) h + L_0 ^*  (f) L (f) h.$$
 We claim that  the equation $ d\tilde f /dt = P (f)\tilde f $ is parabolic   (for a given $f$).  To see this we must
 examine  the symbol
 $$ \sigma P (f) \xi  = \sigma  DE (f) \sigma + \sigma L_0 ^* (f) (\xi) \cdot  \sigma  L_0(f) (\xi).\leqno (6. 13) $$
 Suppose $v$ is an eigenvector  in $F$ with eigenvalue $\lambda$. Then  $\sigma P (f)  (\xi) v = 0$. But
 $\sigma L_0 (f) (\xi) \cdot \sigma  DE (f) (\xi) = 0$, so applying $\sigma L_0 (f)$ to the LHS and RHS of (6.13) we get
 
 $$ \sigma L_0 (f) (\xi) \cdot \sigma L_0  ^* (f) \xi \cdot \sigma L_0 (f) (\xi ) v =  \lambda \sigma L_0 (f) (\xi) v.$$
 Taking inner product of  the   above equality  with $\sigma L_0 (f) (x) v $ we get 
 $$ | \sigma L_0 ^*  (f) (\xi) \cdot \sigma L_0(f) (\xi)  v |  ^2 = \lambda |\sigma L_0(f) (\xi) v | ^2.$$
 
 Now  if $\sigma  L_0 ^* (f) \cdot \sigma L_0 (f) (\xi)  v = 0$ then $\sigma L_0 (f) (\xi ) v = 0$, and otherwise
 $\lambda $ is real and strictly positive. When  $\sigma L_0 (f)  (\xi) v = 0$, then  $\sigma DE (f) (\xi ) v = \lambda v$
 by  (6.13)  and $\lambda $ has  strictly positive real part by our hypothesis in Theorem 6.6.  Thus $P(f)$ is parabolic.
 
 We proceed to solve the system of equations
 $${ d\tilde f \over dt} - P(f) \tilde f+ L_0 ^* (f) \tilde g = \tilde h, \leqno (6.14.1)$$
 $${ d \tilde g \over dt } - M (f)\tilde f = \tilde k\leqno (6.14.2)$$
 for the unknown function $\tilde  f $ and  $\tilde g$ for given $\tilde h$ and $\tilde k$ and
 given $f$, with initial data $ \tilde f = \tilde f_0$ and $\tilde g  =\tilde g_0 = L (f_0)  \tilde f_0$ at $t = 0$.
 
 In Step 3 below  we prove that the solution  $(\tilde f, \tilde g)$ exists  and is unique,
 and is  a smooth tame function of $(f, \tilde h, \tilde k,  \tilde f_0, \tilde g_0)$.
Then putting $\tilde l = \tilde g  - L(f) \tilde f$  and substituting $\tilde k = L_0 (f)   \tilde h$ we get
$$ {d\tilde l \over dt } = {d\tilde  g \over dt} - L_0 (f)  { d\tilde f \over dt }$$
$$ = L _0 (f) DE (f) \tilde f + \tilde k  - L_0 (f)   { d\tilde f \over dt} $$
$$ \stackrel{(6.14.1)}{= } -L_0 (f) DE (f)  \tilde f - L_0 (f) P(f) \tilde f +L_0 (f) L_0^* (f)\tilde g$$
$$ \stackrel{(6.13)}{=} L_0(f) [ - L_0 ^* (f) L(f) \tilde f   + L_0 ^* (f) L(f) (\tilde l  + L (f) \tilde f)] $$
$$ = L_0 (f) L_0 ^* (f) \tilde l,\leqno (6.15)$$
 and $\tilde l   = 0$ at $t = 0$.  But then  (6.15) implies the obvious integral  inequality
 $${d \over dt}\int _X | \tilde l |  ^2  d\mu + 2 \int_X | L_0 ^* (f) \tilde l | d \mu = 0.$$
 Hence  $ \tilde  l = 0$. Then it follows that  $ \tilde g = L(f) \tilde f$.  Using this and  we get from (6.14.1)
 $${d\tilde f \over dt } - DE (f) \tilde f  = \tilde h.$$
 This completes  Step 2.
 
 \medskip
 
 STEP 3. {\it  The system (6.14.1) and (6.14.2) is a weakly parabolic linear system whose  smooth solution  uniquely exists.}
 
 Set $ P_0 (f)h: =  DE (f)h +  L_0 ^* (f) L_0 (f) h.$
 Then $P_0(f)h$ is a linear   differential operator in $h$  and $ P(f) h = P _0 (f) h +  L_0  ^* (f) A (f) $.
 Set $ h = \tilde h - L_0 ^* (f)  A (f)$.  Since $f$  in the system of (6.14.1) and (6.14.2) is given,  we shall re-denote a given constant $\tilde k$ by $k$, 
 variables $\tilde f, \tilde g$,  by $f, g$  and linear differential operators  $P_0 (f), L_0 ^* (f), M(f)$ by
 $ P,  L, M$.  Then the  system of (6.14.1) and (6.14.2) is equivalent to the following system of linear  evolution equations on
 $ 0 \le t \le T$ for sections $f $ of $F$ and $g $ of $G$
 $$ {d f \over dt } = P f + L g  + h, \hspace {1cm} {d g \over dt} = M f +k.\leqno (6.16)$$
   
   Clearly the existence, uniqueness and smoothness of a solution of (6.16) is a consequence  of  Hamilton's theorem
   [Hamilton1982, Theorem 6]. He considered the following equation
   $${ df\over dt } = P f + L g + h, \hspace {1cm} { dg \over dt } = Mf + Ng  + k \leqno (6.17) $$
   where $P, L, M $ and $N$ are linear differential operators involving only space derivatives whose coefficients are
   smooth functions of both  space and time. He assumed that $P$ has degree 2, $L$ and $M$ have degree 1  and $N$ has degree
   0.
   
   \medskip
   
 {\bf 6.18. Theorem} (\cite{Hamilton1982}, Theorem 6). { \it  Suppose the equation $ df / dt  = Pf $ is parabolic.
 Then for any given $(f_0, g_0, h, k)$ there exists  a unique   smooth  solution $ (f, g)$ of the system
 (6.17) with $f = f _0$ and $g = g _0$ at $t = 0$.}
 
 \medskip
 
 The proof of this Theorem   occupies the whole section 6 in   Hamilton's paper.

 \medskip
 
  Finally  we  formulate a conjecture which might be solved by using  the  Yang-Mills bar equation and might be helpful  for understanding the Hodge conjecture.
  A unitary connection $A$ on a Hermitian bundle $E$ over a projective algebraic manifold $M$ is holomorphic, if  the $L ^ p_q $-norm
  of the  component $F^{0,2}_A$ less than some positive constant $\eps(M)$, where $p,q$ are  some integers depending on   the dimension of $M$.  
  
  In  a subsequent paper we shall show the long time existence of  a Yang-Mills bar gradient-like
  flow and    discuss its consequences.

 \medskip

 {\bf Acknowledgement}.  This paper is partially supported by grant of ASCR Nr IAA100190701.  A large part of the paper has been conceived during   my visits
 to the Max-Planck-Institute for Mathematics in Leipzig in last years.  I thank
 J\"urgen Jost for his support and   stimulating discussions.

\bigskip


\end{document}